\documentstyle[12pt,twoside]{article}
\input{amssym.def}
\input{amssym}
\newtheorem{theorem}{Theorem}

\begin{document}
\def\be{\begin{equation}}
\def\ee{\end{equation}}
\def\T{{\Bbb T}}
\def\normL2#1{\|{#1}\|_{L^2}}
\def\normomL2#1{\|{#1}\|_{L^2(\Omega)}}
\def\normdelta#1{\|{#1}\|^2_\delta}
\def\normepsilon#1{\|{#1}\|^2_\epsilon}
\def\normtwodelta#1{\|{#1}\|^2_{2\delta}}
\def\norms#1#2{\|{#1}\|^2_{#2}}
\def\ip#1#2{({#1},{#2})_{L^2}}
\def\boxb{\Box _b}
\title{Global (and Local) Analyticity for Second Order
Operators  Constructed from Rigid Vector Fields on
Products of Tori}
\author{David S. Tartakoff\\ Department
of Mathematics\\University of Illinois at
Chicago\\ 851 S. Morgan St., m/c 249\\
Chicago Illinois 60607-7045, U.S.A.}
\date{}
\maketitle
\begin{abstract}
We prove global analytic hypoellipticity on a product of tori
for partial differential operators which are constructed as
rigid (variable coefficient) quadratic polynomials in real
vector fields satisfying the H\"ormander condition
and where $P$ satisfies a `maximal' estimate.
We also prove an analyticity result that is local in some
variables and global in others for operators whose
prototype is
 $$ P= \left({\partial \over {\partial
x_1}}\right)^2 +  \left({\partial \over {\partial x_2}}\right)^2
+  \left(a(x_1,x_2){\partial \over {\partial t}}\right)^2.$$
(with analytic $a(x), a(0)=0,$ naturally, but not identically
zero). The results, because of the flexibility of the methods,
generalize recent work of Cordaro and Himonas in
\cite{Cordaro-Himonas 1994} and Himonas in
\cite{Himonas 199X} which showed that certain operators
known not to be locally analytic hypoelliptic (those of
Baouendi and Goulaouic
\cite{Baouendi-Goulaouic 1971}, Hanges and Himonas
\cite{Hanges-Himonas 1991}, and Christ \cite{Christ 1991a})
were {\it globally} analytic hypoelliptic on products of tori.
\end{abstract}

\section{Introduction}
We are concerned here with the global
analytic hypoellipticity of second
order operators of the form
   \be P = \sum_{j,k =
   1}^ma_{jk}(x,t)X_jX_k + \sum_{j =
   1}^mb_j(x,t)X_j + b_0(t)X_0 +
c(x,t)\label{1}\ee
on a product of tori,
  \be\T^N = \T^m \times \T^n,\label{2}\ee
where $x \in \T^m, t \in \T^n.$
   Here the functions
$a_{jk},b_j \,$ and $c(x,t)$
may be complex valued, but the
`rigid' vector fields
   \be X_j = \sum_{k=1}^m d_{jk}(x)
   {\partial\over \partial x_k} +
   \sum_{k=1}^n e_{jk}(x){\partial
   \over \partial t_k}\label{3}\ee
are real.
   The hypotheses we make are that 	
   \be\{X_j^\prime = \sum_{k=1}^m
   d_{jk}(x)
   {\partial\over \partial
x_k}\}\label{4}\ee are
independent, $k=1,\ldots,m,$
and there exists a constant C such
that for all smooth $v$,
   \be\sum_{j=0}^m\normL2{X_jv}^2 + \normL2{v}^2
 \leq
   C\{|\Re\ip{Pv}{v}| +
   \norms{v}{-1}\}.\label{5}\ee
For example, if
the vector fields
   $\{X_j\}_{j=0,...,m}$
satisfy the H\"ormander condition
that their iterated
commutators span then whole tangent
space and the matrix
   $A=(a_{jk})$
is the identity, then one even has a
subelliptic estimate, which implies
{\sl arbitrary positivity} (an arbitrarily
large multiple of the second term on
the left, provided one adds a
corresponsingly sufficiently large
multiple of the negative norm on the
right). The positivity of the
self-adjoint matrix $A = (a_{jk})$ alone
will give an estimate of this form without
the second term on the left and with the
norm on the right replaced by the $L^2$
norm, but we need very slightly more. For
example, the positivity of  $A$ together
with
  \be -\sum (X_ja_{jk})X_k-\sum
  ((X_ja_{jk})X_k)^*
  +\sum b_j(x,t)X_j$$
  $$+ (\sum
  b_j(x,t)X_j)^* + X_0 + X_0^*
  + c(x,t) > 0\label{6}\ee
will suffice, and this in turn would
follow from sufficient positivity of the
zero order term $c(x,t).$
\par
This class of operators generalizes that given in
\cite{Cordaro-Himonas 1994}, and in our opinion
simplifies the proof. The more flexible proof techniques we
employed in  \cite{Tartakoff 1976} and \cite{Tartakoff 1978}
allow us to handle this broader class of operators. At one
point in \cite{Cordaro-Himonas 1994} the authors also
prove a theorem of analyticity that is global in some
variables and local in others for operators like
\be P= \left({\partial \over {\partial x_1}}\right)^2 +
\left({\partial \over {\partial x_2}}\right)^2 +
\left(a(x){\partial \over {\partial t}}\right)^2.\ee
Our methods apply to these operators as well (cf.
Theorems 2, 3 and 4).
\par Our interest in these probems was stimulated by the
work of Cordaro and Himonas \cite{Cordaro-Himonas 1994}.
\section{Statement and Proofs of the Theorems}
\begin{theorem} Let $P$ be a partial differential operator of
the form (1) above with real analytic coefficients
$a_{jk}(x,t), b_k(x,t)$ and $c(x,t),$
where the real analytic vector fields
$\{X_j\}_{j=0,\ldots ,n}$ are `rigid' in
the sense of (4). Assume that $P$ satisfies
the {\it a priori} estimate (5) for some
$C\geq 0.$ Then P is globally analytic
hypoelliptic - that is, if $v$ is a
distribution on ${\Bbb T} ^N,$ with $Pv$ analytic
on ${\Bbb T} ^N,$ then $v$ itself is analytic on
${\Bbb T} ^N.$
\end{theorem}

We also state three theorems which are local in some
variables and global in others. In so doing, we hope to
elucidate the distinction between local and global
analyticity. These results are stated for rather explicit, low
dimensional operators for easy reading. For much fuller
and more general results, the reader is referred to the
forthcoming paper of Bove and Tartakoff
(\cite{Bove-Tartakoff 1994}). The restriction to second order
operators is undoubtedly artificial, as the methods of our
recent paper with Popivanov
\cite{Popivanov-Tartakoff 1994} suggest.

First we assume that
$x \in {\Bbb T}
^2,$ but that $t \in I,$ where $I$ an open interval:
\begin{theorem} Let the operator $P$ be given
by
\be P= \left({\partial \over {\partial x_1}}\right)^2 +
\left({\partial \over {\partial x_2}}\right)^2 +
\left(a(x_1,x_2){\partial \over {\partial
t}}\right)^2=\sum_1^3 X_j^2.\ee with $x\in {\Bbb T} ^2$ but
$t\in I, I$ an interval. Then if
$a(x_1,x_2)$ is analytic, zero at $0$ but not identically zero
(so that the H\"ormander condition is satisfied for $P$), and
$Pu=f$ with $f$ real analytic on ${\Bbb T}^2 \times I,$
then $u$ is also analytic on ${\Bbb T} ^2 \times I.$
\end{theorem}
\par\noindent
{\bf Remark.} Theorem 2 holds for a wide class of
operators of this type. For example, if we denote by $Y_j$
the vector fields
$$Y_1 = {\partial \over {\partial x_1}},\,\, Y_2 = {\partial
\over {\partial x_2}}, \hbox{ and } Y_3=a(x_1,x_2){\partial
\over {\partial t}}$$
then Theorem 2 holds for any second order polynomial in the
$Y_j$  \be P = \sum_{|\alpha|\leq 2}b_\alpha (x,t)
Y_{I_\alpha} \ee
with
(non-rigid) variable coefficients $b_\alpha (x,t)$ such that
(5) holds with $X_j$ replaced by $Y_j.$

Next we look at what happens with
$x_1 \in I_1,$ $x_2 \in {\Bbb T}
^2,$ and $t \in I_2,$ when the coefficient $a(x) = a(x_1),$
where the
$I_j$ are open intervals:
\begin{theorem} Let the operator $P$ be given
by
\be P= \left({\partial \over {\partial x_1}}\right)^2 +
\left({\partial \over {\partial x_2}}\right)^2 +
\left(a(x_1){\partial \over {\partial
t}}\right)^2=\sum_1^3 X_j^2\ee
with $x_1 \in I_1,$ $x_2 \in
{\Bbb T} ^1,$ and $t \in I_2,$ the $I_j$ being intervals. Then
if
$a(x_1)$ is analytic, zero at $0$ but not identically zero (so
that the H\"ormander condition is satisfied for $P$), and
$Pu=f$ with $f$ real analytic on $I_1\times{\Bbb T} ^1
\times I_2,$ then $u$ is also analytic on $I_1\times{\Bbb T}
^1 \times I_2.$
\end{theorem}

Finally we consider the case where $a(x) = a(x_1,x_2),$
not identically zero,
has the form
$$a^2(x) = a_1^2(x_1)+a^2_2(x_2),$$
\begin{theorem} Let the operator $P$ be given
by
\be P= \left({\partial \over {\partial x_1}}\right)^2 +
\left({\partial \over {\partial x_2}}\right)^2 +
\left(a_1^2(x_1)+a_2^2(x_2)\right)\left({\partial \over
{\partial t}}\right)^2=\sum_1^4 X_j^2\ee
with $x\in {\Bbb T} ^m$ but $t\in I, I$ an
interval. Then if
$a(x)$ is analytic, zero at $0$ but not identically zero (so
that the H\"ormander condition is satisfied for $P$), and
$Pu=f$ with $f$ real analytic near $0,$ then $u$ is
real analytic near $0.$
\end{theorem}
\par\noindent
{\bf Remark } These theorems have evident microlocal
versions and allow suitable variable coefficient
combinations of the appropriate vector fields as well as the
addition of lower order terms in these vector fields
(${\partial
\over {\partial x_1}}, {\partial \over {\partial x_2}},$
and $a(x){\partial \over {\partial t}}$ for Theorem 2,
$a_1(x_1){\partial
\over {\partial t}}$ in the case of Theorem 3, and
$a_1(x_1){\partial
\over {\partial t}}$ {\it and} $a_2(x_2){\partial
\over {\partial t}}$ in the case of the Theorem 4).

\section{Proofs of the Theorems}
For the moment we shall assume that $v$ and $u$ are
known to belong to $C^\infty ,$ and at the
end make some comments about the $C^\infty$
regularity of the solutions.
\subsection{Proof of Theorem 1}
Using well known results, it suffices to
show that, in $L^2$ norm, we have Cauchy
estimates on derivatives of $v$ of the form
   \be\normL2{X^\alpha T^\beta v} \leq
   C^{|\alpha| + |\beta| +1}|\alpha|!
	  |\beta|!\label{7}\ee
for all $\alpha$ and $\beta$. And
microlocally, since the operator is
elliptic in the complement of the span
$W$ of the vector fields $\partial \over
\partial t_j ,$ it suffices to look near
$W$, and there all derivatives are bounded
by powers of the $\partial \over \partial
t_j ,$ alone. That is, modulo analytic
errors, we may take $\alpha = 0$ and indeed
$\beta = (0,\ldots,0,r,0,\ldots,0),$
as follows by integration by parts and a
simple induction.
Here $T = (T_1,\ldots,T_m)$
with
\be T_k = {\partial \over {\partial
    t_k}},\label{8}\ee		
$k = 1, \ldots, m,$ and we shall take
$T^\beta = T_1^b$ for simplicity.  In
particular, note that
  \be [T_l, P] = \sum_{j,k =
   1}^ma_{jk}^\prime(x,t)X_jX_k + \sum_{j =
   1}^mb_j^\prime(x,t)X_j +
		c^\prime(x,t)\label{9}\ee
and thus
 \be |\Re\ip{[P,T_1^b]v}{T_1^bv}|
  \leq C|\sum_{b^\prime \geq
   1}{b\choose{b^\prime}}
   \ip{c^{(b^\prime)}X^2T_1^{b-b^\prime}v}
   {T^bv}| \label{10}\ee
  $$\leq |\sum_{b^\prime \geq
   1}{b\choose{b^\prime}}
   \ip{c^{(b^\prime
   +1(-1))}XT_1^{b-b^\prime}v}
   {(X)T_1^bv}|,$$
 $$\leq l.c.\sum_{b^\prime \geq1}b^{2b^\prime}
   C_c^{2(b^\prime +1)}
   \normL2{XT_1^{b-b'}v}^2
   + s.c.\normL2{(X)T_1^bv}^2,$$
where the $(X)$ on the right represents the
fact that this X may not or may be present,
depending on whether $c^{(b')}$
received one more derivative or not. We have
written $X^2$ for a generic $X_jX_k$ as well
as $C_c$ for the largest of the constants which
appear in the Cauchy estimates for the analytic
coefficients in $P.$ The large constant ($l.c.$)
and small constant ($s.c.$) are independent of
$b$, of course, the small constant being small
enough to allow this term to be absorbed on
the left hand side of the inequality. Absorbing
yields:

  $$\sum_1^m\normL2{X_jT_1^bu}^2
    + \normL2{X_jT_1^bu}^2
    \leq C\{\normL2{T^bPu}^2
   + \sum_{b\geq b'\geq 1} b^{2b^\prime}
   C_c^{2(b^\prime +1)}
   \normL2{XT_1^{b-b'}v}^2\}$$
which, iterated until the last term on the right
is missing, gives

  $$\sum_1^n\normL2{X_jT_1^bu}^2
    + \normL2{T_1^bu}^2
    \leq C_{(Pu)}^{b+1}b!,$$
which implies the analyticity of $u.$

Finally, we have taken the solution to belong to $C^\infty;$
for the $C^\infty$ behavior of the solution, the methods of
\cite{Tartakoff 1973} which utilize the hypoellipticity
techniques of \cite{Hormander 1963} will lead
quickly to the $C^\infty$ result.

\subsection{Proof of Theorem 2}
We recall that we are (for simplicity) taking the operator $P$
to have the particular form \be P= \left({\partial \over
{\partial x_1}}\right)^2 +  \left({\partial \over {\partial
x_2}}\right)^2 +  \left(a(x){\partial \over {\partial
t}}\right)^2.\ee
Again we take $u \in C^\infty,$ since
subellipticity is a local phenomenon and the operators we
are dealing with are clearly subelliptic under our
hypotheses.  And again it suffices to estimate derivatives in
$L^2$ norm, i.e., to show that with $\phi , \psi$ of compact
support, and
\be X_j = {\partial \over {\partial x_1}},  {\partial \over
{\partial x_2}}, {\hbox{ or }} a(x){\partial \over {\partial
t}},\label{def:X_j}\ee
we will have
\be
\sum_j\normL2{X_j\psi (x)\phi (t) Z^p u} \leq C_u^{p+1}p!
\ee where each $Z$ is (also) of the form
$Z={\partial \over {\partial x_1}},  {\partial \over {\partial
x_2}},$ or  $a(x){\partial \over {\partial t}}.$ That this will
suffice is due to a result by Helffer and Mattera
(\cite{Helffer-Mattera 1980}) but it won't save us any work
as we find (\cite{Tartakoff 1976}, \cite{Tartakoff 1978},
\cite{Tartakoff 1978}) that in trying to bound powers of the
first two types of $Z$ we are led to needing to establish
analytic type growth of derivatives measured by powers of
$Z$ of the form ${{\partial}\over{\partial t}}$ itself.  Actually
we shall show that for any given $N,$  there exists a
localizing function $\phi_N(t) \in C_0^{\infty}$ and $\psi (x)$
(independent of $N$) with
\be \sum_j\normL2{X_j\psi (x)\phi_N(t)
Z^p u} \leq C_u^{p+1}N^p, \qquad p \leq N. \ee And in fact
the functions $\phi_N(t)$ will be chosen to satisfy \be |D^r
\phi_N(t)| \leq C_u^{r+1}N^r, \qquad r\leq N \ee  uniformly
in $N.$
\par The philosophy of all $L^2$ proofs is to replace
$v$ in (5) by  $\psi (x)\phi (t) \tilde{Z}^p u$ with
$\tilde{Z}={\partial \over {\partial x_1}},  {\partial \over
{\partial x_2}},$ or  ${\partial \over {\partial t}}$ and
commute $\psi (x)\phi (t) \tilde{Z}^p$ past the differential
operator $P$. For argument's
sake, and since everything else is simpler, we may restrict
ourselves to the worst case which is given by
$\tilde{Z}={\partial \over {\partial t}}.$ In doing so, we
encounter the errors
\be [{\partial\over {\partial x_1}},\psi
\phi ({\partial \over {\partial t}})^p], \,\,
[{\partial\over {\partial x_2}},\psi\phi({\partial
\over {\partial t}})^p], {\hbox{ and }}
[a(x){\partial\over {\partial t}},\psi\phi({\partial
\over {\partial t}})^p] \ee
Thus, starting with a
given value of $p$, the left hand side of the {\it a priori}
inequality (5) will bound $X_j\psi
(x)\phi (t) ({\partial \over {\partial t}})^pu$ in $L^2$ norm
(after taking the inner product on the right and integrating
by parts one derivative to the right) by

\be \normL2{{{\partial \psi (x)\phi (t)}\over {\partial x}}
({\partial \over {\partial t}})^pu}
{\hbox{ and }} \normL2{a(x){{\partial \psi (x)\phi (t)} \over
{\partial t}} ({\partial \over {\partial t}})^pu}, \ee
(and related terms arising from the integration by parts,
terms which exhibit the same qualitative behavior as these).

So, at the very least, we have bounded
$ \normL2{X_j\psi (x)\phi (t) ({\partial \over {\partial
t}})^pu}$ by
\be \normL2{\psi' (x)\phi (t)
({\partial \over {\partial t}})^pu}
{\hbox{ and }} \normL2{a(x)\psi (x)\phi' (t)({\partial \over
{\partial t}})^pu}. \label{baderrors}\ee

In the first term, we have
lost the `good' derivative ${\partial \over {\partial
x}}$ and seen it appear on the localizing function, but we
cannot iterate this process, since the {\it a priori} estimate
(5) is only truly effective when a `good' derivative is
preserved; when no `good' derivative is preserved, we have
seen often enough that only in the global situation, when
the derivative that appeared on the localizing function can
be absorbed by a {\em constant} by the introduction of a
partition of unity (in that variable, in this case it is the $x,$
or toroidal, variable) can one obtain analyticity (in {\em
that} variable). The second type of term in (\ref{baderrors})
is actually good, since the factor $a(x)$ will
combine with one of the `bad' derivatives ${\partial \over
{\partial t}}$ to give a `good' derivative
$Z=a(x){\partial \over{\partial t}}$ which may be iterated
under (5). That is, modulo terms which lead to global
analyticity in $x,$ we have the iteration schema
\be\sum_j\normL2{X_j\psi (x)\phi (t) ({\partial \over
{\partial t}})^pu}  \rightarrow
C\sum_j\normL2{X_j\psi (x)\phi' (t) ({\partial \over {\partial
t}})^{p-1}u}\ee
which may indeed be iterated. The result of multiple
iterations is, for $\phi(t) =\phi_N(t)$ satisfying the estimates
(16) but the localizations in $x$ merely smooth and
subject to  $\sum_k\psi_k(x) = 1$ (none will ever receive
more than a couple of derivatives)
$$\sum_{j,k}\normL2{X_j\psi_k (x)\phi_N (t) ({\partial
\over {\partial t}})^pu}\leq $$
$$\leq \sum_{{j,k}\atop{p'\leq p}}
C^{p'}{\normL2{\psi_k (x)\phi_N^{(p-p')}(t)({\partial
\over {\partial t}})^{p'}Pu} +
\sum_{{j,k}\atop{p'\leq p}}
C^pN^{p'}\normL2{(X_j)\psi_k (x)\phi_N^{(p')}(t)u}}.$$

Since $p\leq N$ and $N^N \leq C^{N+1}N!$ by Stirling's
formula, under the bounds (16) this yields the desired
analyticity,  which is local in $t$ but global in
$x.$

\subsection{Proof of Theorem 3}
The new ingredient in Theorem 3 is that the function 	
$a(x)$ is now of a more special form. Thus it is only 	
on the hypersurface $x_1=0$ that the operator $P$ is not
elliptic; if in the above proof we replace the compactly
supported function $\psi (x)$ by a product:
$$\psi (x) = \psi_1(x_1)\psi_2(x_2),$$
with both $\psi_j(s)$ equal to one near $s=0,$ then when
derivatives enter on $\psi_1(x_1),$ the support of
$\psi_1^\prime $ is contained in the elliptic region, and 	
only in $x_2$ does one need to pass to further and further
patches, ultimately using a (finite) partition of unity on the
torus in $x_2.$

\subsection{Proof of Theorem 4}
The new ingredient in Theorem 4 is that there are four
vector fields,
$${\partial \over{\partial x_1}}, {\partial \over{\partial
x_2}}, a_1(x_1){\partial \over{\partial t}},
{\hbox{ and }} a_1(x_1){\partial
\over{\partial t}}.$$
The above considerations apply to $x_1$ and $x_2$
separately, now, since if {\em either} $x_1\neq 0$ {\em or}
$x_2\neq 0$ we are in the elliptic region where the solution
is known to be analytic.

\par\noindent{\bf Remark} It is not hard to see that
derivatives in $x_1$ and $x_2$ {\em always} behave well -
i.e. that $(x_1, x_2, t; \xi_1, \xi_2, \tau)$ is never in the
analytic wave front set $WF_A(u)$ for $(\xi_1, \xi_2) \neq
(0,0)$ whenever this is true (punctually) of $Pu,$ since only
points of the form $(x_1, x_2, t; 0,0,\tau)$ are characteristic
for $P.$ Hence the above theorems are actually
"microlocal(-global)" in a sense which is fairly evident,
much as in \cite{Derridj-Tartakoff 1993c}.


\begin{thebibliography}{99}

\bibitem{Baouendi-Goulaouic 1971}
{\sc M.S. Baouendi and C. Goulaouic}, \
{\it Analyticity for Degenerate
Elliptic Equations and Applications.}
Proc. Symp. in Pure Math., {\bf 23}
(1971), 79-84.

\bibitem{Bove-Tartakoff 1994} {\sc A.
Bove and D. S. Tartakoff}, \
{\it Microlocal Gevrey
Hypoellipticity for Subelliptic
Operators} To Appear.

\bibitem{Christ 1991a} {\sc M. Christ},
\ {\it Certain Sums of Squares of
Vector Fields Fail to Be Analytic
Hypoelliptic} Comm. in P.D.E..,
{\bf 10} no. 10 (1991), 1695-1707.

\bibitem{Cordaro-Himonas 1994} {\sc
P. Cordaro and A. Himonas}, \ {\it
Global Analytic Hypoellipticity of a
Class of Degenerate Elliptic Operators
on the Torus}, Mathematical Research Letters
{\bf 1} (1994), 501-510.

\bibitem{Derridj-Tartakoff 1993c}
{\sc M. Derridj and D.S. Tartakoff},
\ {\it Global Analyticity for $\boxb$
on Three Dimensional Pseudoconvex CR
Manifolds.} Communications in P. D. E.
{\bf 18}(11) 1993, 1847-1868.

\bibitem{Hanges-Himonas 1991}
{\sc N. Hanges and A. Himonas}, \ {\it
Singular Solutions for Sums of
Squares of Vector Fields.} Comm. in
P.D.E. {\bf 16} no. 8 \& 9 (1991),
1503-1511.

\bibitem{Helffer-Mattera 1980} {\sc B. Helffer and C. Mattera},
\ {\it Analyticit\'{e} et
It\'{e}r\'{e}s R\'{e}duits d'un Syst\`{e}me de Champs de Vecteurs.}
Comm. in P.D.E.'s, {\bf 5}
(1980), 1065-1072.

\bibitem{Himonas 199X} {\sc A. Alexandrou Himonas}, \
{On Degenerate Elliptic Operators of Infinite Type}, to
appear, Mathematische Zeitschrift.

\bibitem{Hormander 1963} {\sc L.
Hormander}, \ Linear Partial
Differential Operators, Springer
Verlag, New York 1969.

\bibitem{Popivanov-Tartakoff 1994}
{\sc P. Popivanov and D. S. Tartakoff}
\ {\it Gevrey Hypoellipticity for
Fourth Order Differential Operators.}
To appear, Comm. in P.D.E.

\bibitem{Tartakoff 1973}
{\sc D. S. Tartakoff}, \ {\it Gevrey
Hypoellipticity for Subelliptic
Boundary Value Problems.} Comm. on Pure
Appl. Math., {\bf 26} (1973), 251-312.

\bibitem{Tartakoff 1976}
{\sc D. S. Tartakoff}, \ {\it On the
Global Real Analyticity of Solutions to
$\boxb$ on Compact Manifolds.} Comm. in
P.D.E., {\bf 1} (1976), 283-311.

\bibitem{Tartakoff 1978} {\sc D. S. Tartakoff},
\ {\it Local Analytic
Hypoellipticity for $\boxb$ on Non-Degenerate
Cauchy Riemann Manifolds.} Proc. Nat. Acad. Sci.
U.S.A. {\bf 75} no. 7 (1978), 3027-3028.

\end{thebibliography}
\end{document}